\documentclass[12pt,thmsa]{article}%
\usepackage{amsmath}
\usepackage{amssymb}
\usepackage{sw20lart}
\usepackage{amsfonts}
\usepackage{graphicx}%
\begin{document}

\title{Metric and Gauge Extensors}
\author{{\footnotesize A. M. Moya}$^{2}${\footnotesize , V. V. Fern\'{a}ndez}$^{1}%
${\footnotesize and W. A. Rodrigues Jr}.$^{1}${\footnotesize \ }\\$^{1}\hspace{-0.1cm}${\footnotesize Institute of Mathematics, Statistics and
Scientific Computation}\\{\footnotesize \ IMECC-UNICAMP CP 6065}\\{\footnotesize \ 13083-859 Campinas, SP, Brazil }\\{\footnotesize e-mail: walrod@ime.unicamp.br }\\{\footnotesize \ }$^{2}${\footnotesize Department of Mathematics, University
of Antofagasta, Antofagasta, Chile} \\{\footnotesize e-mail: mmoya@uantof.cl}}
\maketitle

\begin{abstract}
In this paper, the second in a series of eight we continue our development of
the basic tools of the multivector and extensor calculus which are used in our
formulation of the differential geometry of smooth manifolds of arbitrary
topology . We introduce metric and gauge extensors, pseudo-orthogonal metric
extensors, gauge bases, tetrad bases and prove the remarkable \emph{golden}
formula, which permit us to view any Clifford algebra $\mathcal{C}\ell(V,G)$
as a deformation of the euclidean Clifford algebra $\mathcal{C}\ell(V,G_{E})$
discussed in the first paper of the series and to easily perform calculations
in $\mathcal{C}\ell(V,G)$ using $\mathcal{C}\ell(V,G_{E})$.

\end{abstract}
\tableofcontents

\section{Introduction}

This the second paper in a series of eight. As emphasized in \cite{1} the
euclidean geometric algebra plays a key role in our presentation. Here we
continue our formulation of the theory of geometrical algebras and extensors,
introducing in Section 2 the concept of \textit{metric extensor} and
\textit{metric adjoint} operators of arbitrary signatures. In Section 3 we
introduce orthogonal metric extensors and in particular the important case of
Lorentz extensors. Gauge extensors, a fundamental tool in the formulation of
geometric theories of the gravitational field are studied in Section 4. In
Section 5 we exhibit some applications of the formalism, such as \textit{gauge
bases}, \textit{tetrad bases} and some algebraic aspects of the tetrad
formalism. Section 6 introduces an prove the remarkable \emph{golden formula},
which permit us to interpret any arbitrary metric Clifford algebra
$\mathcal{C}\ell(V,G)$ as a deformation of the euclidean algebra
$\mathcal{C}\ell(V,G_{E})$ (see \cite{1}) and to perform calculations in
$\mathcal{C}\ell(V,G)$ using $\mathcal{C}\ell(V,G_{E})$. In section 7 we
present our conclusions.

\section{Metric Extensor}

Let $V$ be a real vector space, $\dim V=n$ and let us endow $V$ with an
arbitrary (but fixed once and for all) euclidean metric $G_{E}.$ Whenever $V$
is equipped with another metric $G$ (besides $G_{E}$) there is an unique
linear mapping (i.e., a $(1,1)$-extensor\footnote{ We follows the nomenclature
of \cite{2}, which present all material on extensors needed for the present
series.}) $g:V\rightarrow V$ such that for all $v,w\in V$
\begin{equation}
v\underset{G}{\cdot}w=g(v)\underset{G_{E}}{\cdot}w,\label{ME.1}%
\end{equation}
where $\underset{G_{E}}{\cdot}$ and $\underset{G}{\cdot}$ denote respectively
the scalar products associated to the euclidean metric structure $(V,G_{E})$
and the metric structure $(V,G).$

Such $g$ is symmetric, i.e., $g=g^{\dagger}$ (the adjoint operator $\left.
{}\right.  ^{\dagger}$ is taken with respect to $(V,G_{E})$)$,$ and
non-degenerate\footnote{For any $t\in ext_{1}^{1}(V):\det[t]\in\mathbb{R}$ is
the unique scalar which satisfies $\underline{t}(I)=\det[t]I$ for all non-zero
pseudoscalar $I\in\bigwedge^{n}V.$}, i.e., $\det[g]\neq0.$ It will be called
the\emph{\ metric extensor} for $G$ (of course, relative to $G_{E}$).

In what follows the euclidean scalar product $\underset{G_{E}}{\cdot}$ and the
scalar product $\underset{G}{\cdot}$ will be denoted by the more convenient
notations: $\underset{}{\cdot}$ (only a dot) and $\underset{g}{\cdot}$ (letter
$g$ under a dot), respectively. Then, the contracted products and the Clifford
product of multivectors $X,Y\in\bigwedge V$ based upon $(V,G_{E})$ will be
denoted by $X\lrcorner Y,$ $X\llcorner Y$ and $XY.$ And, those products based
upon $(V,G)$ will be denoted by $X\underset{g}{\lrcorner}Y,$ $X\underset
{g}{\llcorner}Y$ and $X\underset{g}{}Y.$

The relationship between the euclidean metric algebraic structures $(\bigwedge
V,\cdot)$ and $(\bigwedge V,\lrcorner,\llcorner),$ and the metric algebraic
structures $(\bigwedge V,\underset{g}{\cdot})$ and $(\bigwedge V,\underset
{g}{\lrcorner},\underset{g}{\llcorner})$ is given by the following noticeable formulas.

For any $X,Y\in\bigwedge V$
\begin{equation}
X\underset{g}{\cdot}Y=\underline{g}(X)\cdot Y, \label{ME.2}%
\end{equation}
where $\underline{g}$ is the extended\footnote{Recall that $\underline{t}\in
ext(V)$ is the so-called extended of $t\in ext_{1}^{1}(V),$ i.e.,
$\underline{t}(\alpha)=\alpha,$ $\underline{t}(v)=t(v)$ and $\underline
{t}(X\wedge Y)=\underline{t}(X)\wedge\underline{t}(Y)$ for all $\alpha
\in\mathbb{R},$ $v\in V$ and $X,Y\in\bigwedge V.$} of $g. $

For any $X,Y\in\bigwedge V$
\begin{align}
\underset{g}{X\lrcorner}Y  &  =\underline{g}(X)\lrcorner Y\label{ME.3a}\\
X\underset{g}{\llcorner}Y  &  =X\llcorner\underline{g}(Y). \label{ME.3b}%
\end{align}

\subsection{Euclidean and Metric Adjoint Operators}

Let $\bigwedge_{1}^{\diamond}V$ and $\bigwedge_{2}^{\diamond}V$ be two
subspaces of $\bigwedge V$ such that each of them is any sum of homogeneous
subspaces of $\bigwedge V.$ Let us take any $t\in1$-$ext(\bigwedge
_{1}^{\diamond}V;\bigwedge_{2}^{\diamond}V),$ i.e., $t$ is a linear mapping
from $\bigwedge_{1}^{\diamond}V$ to $\bigwedge_{2}^{\diamond}V$, called a
general $1$-extensor over $V$.\footnote{We are using the nomenclature of
\cite{2}.}

Let us denote respectively by $\left.  {}\right.  ^{\dagger}$ and $\left.
{}\right.  ^{\dagger(g)}$ the adjoint operators taken with respect to
$(V,G_{E})$ and $(V,G).$

The euclidean adjoint \cite{2} of $t,$ namely $t^{\dagger},$ satisfies the
$G_{E}$-scalar product condition
\begin{equation}
t(X)\cdot Y=X\cdot t^{\dagger}(Y).\label{ME.4}%
\end{equation}
The metric adjoint of $t$, namely $t^{\dagger(g)}$, satisfy the analogous
$g$-scalar product condition, i.e.,
\begin{equation}
t(X)\underset{g}{\cdot}Y=X\underset{g}{\cdot}t^{\dagger(g)}(Y).\label{ME.5}%
\end{equation}

We find next the relationship between $t^{\dagger}$ and $t^{\dagger(g)}.$

Let us take $X\in\bigwedge_{1}^{\diamond}V$ and $Y\in\bigwedge_{2}^{\diamond
}V.$ We recall the following properties: the adjoint of an extended equals the
extended of the adjoint, and the adjoint of a composition equals the
composition of the adjoints in reversed order. Thus, by using Eq.(\ref{ME.5}),
Eq.(\ref{ME.2}), the euclidean symmetry of $g,$ i.e., $g=g^{\dagger},$ and
Eq.(\ref{ME.4}) we have that
\begin{align*}
X\underset{g}{\cdot}t^{\dagger(g)}(Y)=t(X)\underset{g}{\cdot}Y  &
\Rightarrow\underline{g}(X)\cdot t^{\dagger(g)}(Y)=\underline{g}\circ
t(X)\cdot Y\\
&  \Rightarrow X\cdot\underline{g}\circ t^{\dagger(g)}(Y)=X\cdot t^{\dagger
}\circ\underline{g}(Y).
\end{align*}
Hence, by non-degeneracy of the euclidean scalar product and recalling that
the inverse of an extended equals the extended of the inverse, we finally get
\begin{equation}
t^{\dagger(g)}=\underline{g}^{-1}\circ t^{\dagger}\circ\underline{g}.
\label{ME.6}%
\end{equation}

\section{Orthogonal Metric Extensor}

Let $\{b_{j}\}$ be any orthonormal basis for $V$ with respect to $(V,G_{E}),$
i.e., $b_{j}\cdot b_{k}=\delta_{jk}.$ Once the Clifford algebra
$\mathcal{C\ell}(V,G_{E})$ has been given we are able to construct exactly $n$
euclidean orthogonal metric extensors with signature $(1,n-1).$ The euclidean
orthonormal eigenvectors \footnote{As well-known, the eigenvalues of any
orthogonal symmetric operator are $\pm1$.} for each of them are just being the
basis vectors $b_{1},\ldots,b_{n}.$

Associated to $\{b_{j}\}$ we introduce the $(1,1)$-extensors $\eta_{b_{1}%
},\ldots,\eta_{b_{n}}$ defined by
\begin{equation}
\eta_{b_{j}}(v)=b_{j}vb_{j}, \label{OM.1}%
\end{equation}
for each $j=1,\ldots,n.$

They obviously satisfy
\begin{equation}
\eta_{b_{j}}(b_{k})=\left\{
\begin{array}
[c]{ll}%
b_{k,} & k=j\\
-b_{k}, & k\neq j
\end{array}
\right.  . \label{OM.2}%
\end{equation}
It means that $b_{j}$ is an eigenvector of $\eta_{b_{j}}$ with the eigenvalue
$+1,$ and the $n-1$ basis vectors $b_{1},\ldots,b_{j-1},b_{j+1,}\ldots,b_{n}$
all are eigenvectors of $\eta_{b_{j}}$ with the same eigenvalue $-1.$

As we can easily see, any two of these $(1,1)$-extensors commutates, i.e.,
\begin{equation}
\eta_{b_{j}}\circ\eta_{b_{k}}=\eta_{b_{k}}\circ\eta_{b_{j}},\text{ for }j\neq
k. \label{OM.3}%
\end{equation}
Moreover, they are symmetric and non-degenerate, and euclidean orthogonal,
i.e.,
\begin{align}
\eta_{b_{j}}^{\dagger}  &  =\eta_{b_{j}}\label{OM.4a}\\
\det[\eta_{b_{j}}]  &  =(-1)^{n-1}\label{OM.4b}\\
\eta_{b_{j}}^{\dagger}  &  =\eta_{b_{j}}^{-1}. \label{OM.4c}%
\end{align}
Therefore, they all are orthogonal metric extensors with signature $(1,n-1).$

The extended of $\eta_{b_{j}}$ is given by
\begin{equation}
\underline{\eta_{b_{j}}}(X)=b_{j}Xb_{j}, \label{OM.5}%
\end{equation}

We can now construct an \emph{euclidean} orthogonal metric operator with
signature $(p,n-p)$ and whose euclidean orthonormal eigenvectors are just the
basis vectors $b_{1},\ldots,b_{n}$. It is defined by
\begin{equation}
\eta_{b}=(-1)^{p+1}\eta_{b_{1}}\circ\cdots\circ\eta_{b_{p}}, \label{OM.6}%
\end{equation}
i.e.,
\begin{equation}
\eta_{b}(a)=(-1)^{p+1}b_{1}\ldots b_{p}ab_{p}\ldots b_{1}. \label{OM.7}%
\end{equation}

It is easy to verify that
\begin{equation}
\eta_{b}(b_{k})=\left\{
\begin{array}
[c]{ll}%
b_{k}, & k=1,\ldots,p\\
-b_{k}, & k=p+1,\ldots,n
\end{array}
\right.  , \label{OM.8}%
\end{equation}
which means that $b_{1},\ldots,b_{p}$ are eigenvectors of $\eta_{b}$ with the
same eigenvalue $+1,$ and $b_{p+1},\ldots,b_{n}$ are eigenvectors of $\eta
_{b}$ with the same eigenvalue $-1.$

It is symmetric and non-degenerate, and orthogonal, i.e.,
\begin{align}
\eta_{b}^{\dagger}  &  =\eta_{b}\label{OM.9a}\\
\det[\eta_{b}]  &  =(-1)^{n-p}\label{OM.9b}\\
\eta_{b}^{\dagger}  &  =\eta_{b}^{-1}. \label{OM.9c}%
\end{align}
So, $\eta_{b}$ is an orthogonal metric extensor with signature $(p,n-p)$.

The extended of $\eta_{b}$ is obviously given by
\begin{equation}
\underline{\eta_{b}}(X)=(-1)^{p+1}b_{1}\ldots b_{p}Xb_{p}\ldots b_{1}.
\label{OM.10}%
\end{equation}

What is the most general orthogonal metric extensor with signature $(p,n-p)$?

To find the answer, let $\eta$ be any orthogonal metric extensor with
signature $(p,n-p).$ The symmetry of $\eta$ implies the existence of exactly
$n$ euclidean orthonormal eigenvectors $u_{1},\ldots,u_{n}$ for $\eta$ which
form just a basis for $V.$ Since $\eta$ is orthogonal and its signature is
$(p,n-p)$, it follows that the eigenvalues of $\eta$ are equal $\pm1$ and the
eigenvalues equation for $\eta$ can be written (re-ordering $u_{1}%
,\ldots,u_{n}$ if was necessary) as
\[
\eta(u_{k})=\left\{
\begin{array}
[c]{ll}%
u_{k}, & k=1,\ldots,p\\
-u_{k}, & k=p+1,\ldots,n
\end{array}
\right.  .
\]

Now, due to the orthonormality of both $\{b_{k}\}$ and $\{u_{k}\},$ there must
be an orthogonal operator\footnote{Recall that an operator on $V$ is just a
$(1,1)$-extensor over $V.$} $\Theta$ such that $\Theta(b_{k})=u_{k},$ for each
$k=1,\ldots,n,$ i.e., for all $a\in V:\Theta(a)=\overset{n}{\underset{j=1}{%
{\displaystyle\sum}
}}(a\cdot b_{j})u_{j}.$

Then, we can write
\begin{align*}
\Theta\circ\eta_{b}\circ\Theta^{\dagger}(u_{k})  &  =\Theta\circ\eta_{b}%
(b_{k})=\Theta(\left\{
\begin{array}
[c]{ll}%
b_{k}, & k=1,\ldots,p\\
-b_{k}, & k=p+1,\ldots,n
\end{array}
\right.  )\\
&  =\left\{
\begin{array}
[c]{ll}%
u_{k}, & k=1,\ldots,p\\
-u_{k}, & k=p+1,\ldots,n
\end{array}
\right.  =\eta(u_{k}),
\end{align*}
for each $k=1,\ldots,n.$ Thus, we have
\begin{equation}
\eta=\Theta\circ\eta_{b}\circ\Theta^{\dagger}. \label{OM.11}%
\end{equation}
By putting Eq.(\ref{OM.6}) into Eq.(\ref{OM.11}) we get
\begin{equation}
\eta=(-1)^{p+1}\eta_{1}\circ\cdots\circ\eta_{p}, \label{OM.12}%
\end{equation}
where each of $\eta_{j}\equiv\Theta\circ\eta_{b_{j}}\circ\Theta^{\dagger}$ is
an euclidean orthogonal metric extensor with signature $(1,n-p).$

But, by using the vector identity $abc=(a\cdot b)c-(a\cdot c)b+(b\cdot
c)a+a\wedge b\wedge c,$ with $a,b,c\in V,$ we can prove that
\begin{equation}
\eta_{j}(v)=\Theta(b_{j})v\Theta(b_{j}), \label{OM.13}%
\end{equation}
for each $j=1,\ldots,p.$

Now, by using Eq.(\ref{OM.13}) we can write Eq.(\ref{OM.12}) in the remarkable
form
\begin{equation}
\eta(v)=(-1)^{p+1}\underline{\Theta}(b_{1}\ldots b_{p})v\underline{\Theta
}(b_{p}\ldots b_{1}).\label{OM.14}%
\end{equation}
Such a pseudo orthogonal metric extensor $\eta$ with the same signature as $g$
is called a \emph{Minkowski extensor}. 

\subsection{Lorentz Extensor}

A $(1,1)$-extensor over $V,$ namely $\Lambda,$ is said to be $\eta$-orthogonal
if and only if for all $v,w\in V$
\begin{equation}
\Lambda(v)\underset{\eta}{\cdot}\Lambda(w)=v\underset{\eta}{\cdot}w.
\label{OE.1}%
\end{equation}

By using Eq.(\ref{ME.5}) and recalling the non-degeneracy of the $\eta$-scalar
product, Eq.(\ref{OE.1}) can also be written as
\begin{equation}
\Lambda^{\dagger(\eta)}=\Lambda^{-1}. \label{OE.2}%
\end{equation}
Or, by taking into account Eq.(\ref{ME.6}), we can still write
\begin{equation}
\Lambda^{\dagger}\circ\eta\circ\Lambda=\eta. \label{OE.3}%
\end{equation}

We emphasize that the $\eta$-scalar product condition given by Eq.(\ref{OE.1})
is logically equivalent to each of Eq.(\ref{OE.2}) and Eq.(\ref{OE.3}).

Sometimes, such a $\eta$-orthogonal $(1,1)$-extensor $\Lambda$ will be called
a \emph{Lorentz extensor }(of course, associated to $\eta$)$.$

\section{Gauge Extensor}

\begin{theorem}
Let $g$  and $\eta$ be a pseudo-orthogonal metric extensors, of the same
signature $(p,n-p).$ Then, there exists a non-singular $(1,1)$-extensor $h$
such that
\begin{equation}
g=h^{\dagger}\circ\eta\circ h.\label{GE.1}%
\end{equation}

Such $h$ is given by
\begin{equation}
h=d_{\sigma}\circ d_{\sqrt{\left\vert \lambda\right\vert }}\circ\Theta
_{uv},\label{GE.3}%
\end{equation}
where $d_{\sigma}$ is a pseudo-orthogonal metric extensor, $d_{\sqrt
{\left\vert \lambda\right\vert }}$ is a metric extensor, and $\Theta_{uv}$ is
a \ pseudo-orthogonal operator which are defined by
\begin{align}
d_{\sigma}(a) &  =\overset{n}{\underset{j=1}{\sum}}\sigma_{j}(a\cdot
u_{j})u_{j}\label{GE.3a}\\
d_{\sqrt{\left\vert \lambda\right\vert }}(a) &  =\overset{n}{\underset
{j=1}{\sum}}\sqrt{\left\vert \lambda_{j}\right\vert }(a\cdot u_{j}%
)u_{j}\label{GE.3b}\\
\Theta_{uv}(a) &  =\overset{n}{\underset{j=1}{\sum}}(a\cdot u_{j}%
)v_{j},\label{GE.3c}%
\end{align}
where $\sigma_{1},\ldots,\sigma_{n}$ are real numbers with $\sigma_{1}%
^{2}=\cdots=\sigma_{n}^{2}=1,$ $\lambda_{1},\ldots,\lambda_{n}$ are the
eigenvalues of $g,$ and $u_{1},\ldots,u_{n}$ and $v_{1},\ldots,v_{n}$ are
respectively the orthonormal eigenvectors of $\eta$ and $g.$
\end{theorem}

\textbf{Proof. }As we can see, $d_{\sigma}=d_{\sigma}{}^{\ast},$ and
$d_{\sigma}=d_{\sigma}{}^{\dagger}$ and $\det[d_{\sigma}]=\sigma_{1}%
\ldots\sigma_{n}\neq0,$ thus $d_{\sigma}$ is a pseudo-orthogonal metric
extensor. $d_{\sqrt{\left\vert \lambda\right\vert }}$ is a metric extensor
since $d_{\sqrt{\left\vert \lambda\right\vert }}$ is symmetric, i.e.,
$d_{\sqrt{\left\vert \lambda\right\vert }}=d_{\sqrt{\left\vert \lambda
\right\vert }}^{\dagger}$, and non-degenerate, i.e., $\det[d_{\sqrt{\left\vert
\lambda\right\vert }}]=\sqrt{\left\vert \lambda_{1}\ldots\lambda
_{n}\right\vert }\neq0$. And, $\Theta_{uv}$ is an pseudo-orthogonal operator,
i.e., $\Theta_{uv}=\Theta_{uv}^{\ast},$ since it is the changing basis
extensor between the orthonormal bases $\{u_{k}\}$ and $\{v_{k}\},$ i.e.,
$\Theta_{uv}(u_{k})=v_{k}.$

The non-singularity of $h$ can now be easily proved. We have indeed that
$\det[h]=\det[d_{\sigma}]\det[d_{\sqrt{\left|  \lambda\right|  }}]\det
[\Theta_{uv}^{\dagger}]=\sigma_{1}\ldots\sigma_{n}\sqrt{\left|  \lambda
_{1}\ldots\lambda_{n}\right|  }(\pm1)\neq0.$

Now, by using $\Theta_{uv}(u_{k})=v_{k}$ and the eigenvalue equations for all
of the metric extensors:

$d_{\sqrt{\left\vert \lambda\right\vert }}(u_{k})=\sqrt{\left\vert \lambda
_{k}\right\vert }u_{k},$ $\eta(u_{k})=\left\{
\begin{array}
[c]{ll}%
u_{k}, & k=1,\ldots,p\\
-u_{k}, & k=p+1,\ldots,n
\end{array}
\right.  ,$ $d_{\sigma}(u_{k})=\sigma_{k}u_{k},$ and $g(v_{k})=\lambda
_{k}v_{k}=\left\{
\begin{array}
[c]{ll}%
\left\vert \lambda_{k}\right\vert v_{k}, & k=1,\ldots,p\\
-\left\vert \lambda_{k}\right\vert v_{k}, & k=p+1,\ldots,n
\end{array}
\right.  $, a straightforward calculation yields
\begin{align*}
h^{\dagger}\circ\eta\circ h(v_{k})  &  =\Theta_{uv}\circ d_{\sqrt{\left\vert
\lambda\right\vert }}\circ d_{\sigma}\circ\eta\circ d_{\sigma}\circ
d_{\sqrt{\left\vert \lambda\right\vert }}\circ\Theta_{uv}^{\dagger}(v_{k})\\
&  =\sigma_{k}\sqrt{\left\vert \lambda_{k}\right\vert }\Theta_{uv}\circ
d_{\sqrt{\left\vert \lambda\right\vert }}\circ d_{\sigma}(\left\{
\begin{array}
[c]{ll}%
u_{k}, & k=1,\ldots,p\\
-u_{k}, & k=p+1,\ldots,n
\end{array}
\right.  )\\
&  =\sigma_{k}^{2}(\sqrt{\left\vert \lambda_{k}\right\vert })^{2}\left\{
\begin{array}
[c]{ll}%
v_{k}, & k=1,\ldots,p\\
-v_{k}, & k=p+1,\ldots,n
\end{array}
\right.  =g(v_{k}).
\end{align*}
Hence, since $v_{1},\ldots,v_{n}$ are basis vectors for $V,$ $h^{\dagger}%
\circ\eta\circ h=g.\blacksquare$

It should be noticed that such a $(1,1)$-extensor $h$ satisfying
Eq.(\ref{GE.1}) is not unique. If there is some $h$ satisfying Eq.(\ref{GE.1}%
), then $h^{\prime}\equiv\Lambda\circ h,$ where $\Lambda$ is a $\eta
$-orthogonal extensor, i.e., $\Lambda^{\dagger}\circ\eta\circ\Lambda=\eta$,
also satisfies Eq.(\ref{GE.1}).

Given any metric extensor $g$, any non-singular $(1,1)$-extensor $h$ which
satisfies $g=h^{\dagger}\circ\eta\circ h,$ for some pseudo-orthogonal metric
extensor $\eta$ with the same signature as $g,$ will be called a \emph{gauge
extensor} for $g$.

\begin{remark}
By comparing \emph{Eq.(\ref{GE.1})} with \emph{Eq.(\ref{OE.3})} we can easily
see that a $\eta$-orthogonal extensor $\Lambda$ is just being a gauge extensor
for $\eta$ itself.

In this particular case, \emph{Theorem 1} is reduced to $\eta=d_{\sigma
}^{\dagger}\circ\eta\circ d_{\sigma},$ where $d_{\sigma}$ is the orthogonal
metric extensor given by \emph{Eq.(\ref{GE.3a})}.
\end{remark}

Theorem 1 implies that given an orthogonal metric extensor $\eta$ with
signature $(p,n-p)$ we can indeed construct a metric extensor $g$ which is not
necessarily orthogonal but has the same signature $(p,n-p)$.

Let us take $h\equiv d_{\rho}\circ\Phi,$ where $d_{\rho}$ is the metric
extensor defined by $d_{\rho}(a)=\overset{n}{\underset{j=1}{%
{\displaystyle\sum}
}}\rho_{j}(a\cdot u_{j})u_{j},$ with all of the real numbers $\rho_{1}%
,\ldots,\rho_{n}\neq0$ and $u_{1},\ldots,u_{n}$ being the orthonormal
eigenvectors of $\eta$, and $\Phi$ is an orthogonal operator, i.e.,
$\Phi^{\dagger}=\Phi^{-1}.$

Then, $g\equiv h^{\dagger}\circ\eta\circ h$ is a metric extensor which has the
following properties: the $p$ positive real numbers $\rho_{1}^{2},\ldots
,\rho_{p}^{2}$ are the eigenvalues of $g$ with the associated orthonormal
eigenvectors $\Phi^{\dagger}(u_{1}),\ldots,\Phi^{\dagger}(u_{p}),$ and the
$n-p$ negative real numbers $-\rho_{p+1}^{2},\ldots,-\rho_{n}^{2}$ are the
eigenvalues of $g$ with the associated orthonormal eigenvectors $\Phi
^{\dagger}(u_{p+1}),\ldots,\Phi^{\dagger}(u_{n}).$ Hence, we see that $g$ has
signature $(p,n-p).$

If we want that $g$ has the pre-assigned eigenvalues $\lambda_{1}%
,\ldots,\lambda_{n}$ with the associated orthonormal eigenvectors
$v_{1},\ldots,v_{n},$ then we should choose $\rho_{k}=\pm\sqrt{\left|
\lambda_{k}\right|  },$ and $\Phi$ as defined by $\Phi(a)=\overset
{n}{\underset{j=1}{%
{\displaystyle\sum}
}}(a\cdot v_{j})u_{j}.$

\section{Some Applications}

Let $\{e_{k}\}$ be any basis for $V,$ and $\{e^{k}\}$ be its euclidean
reciprocal basis for $V$, i.e., $e_{k}\cdot e^{l}=\delta_{k}^{l}.$ Let us take
a non-singular $(1,1)$-extensor $\lambda.$ Then, it is easily seen that the
$n$ vectors $\lambda(e_{1}),\ldots,\lambda(e_{n})\in V$ and the $n$
vectors\footnote{Recall that $\lambda^{*}=(\lambda^{-1})^{\dagger}%
=(\lambda^{\dagger})^{-1}.$} $\lambda^{\ast}(e^{1}),\ldots,\lambda^{\ast
}(e^{n})\in V$ define two well-defined euclidean reciprocal bases for $V,$
i.e.,
\begin{equation}
\lambda(e_{k})\cdot\lambda^{\ast}(e^{l})=\delta_{k}^{l}.\label{RT.1}%
\end{equation}
The bases $\{\lambda(e_{k})\}$ and $\{\lambda^{\ast}(e^{k})\}$ are
conveniently said to be a $\lambda$\emph{-deformation} of the bases
$\{e_{k}\}$ and $\{e^{k}\}.$ Sometimes, the first ones are named as the
$\lambda$\emph{-deformed bases }of the second ones.

\subsection{Gauge Bases}

Let $h$ be a gauge extensor for $g,$ and $\eta$ be a pseudo-orthogonal metric
extensor with the same signature as $g.$ According to Eq.(\ref{GE.1}) the
$g$-scalar product and $g^{-1}$-scalar product are related to the $\eta
$-scalar product by the following formulas
\begin{align}
X\underset{g}{\cdot}Y &  =\underline{h}(X)\underset{\eta}{\cdot}\underline
{h}(Y),\label{GB.1a}\\
X\underset{g^{-1}}{\cdot}Y &  =\underline{h}^{\ast}(X)\underset{\eta}{\cdot
}\underline{h}^{\ast}(Y).\label{GB.1b}%
\end{align}

The $\eta$-deformed bases $\{h(e_{k})\}$ and $\{h^{\ast}(e^{k})\}$ satisfy the
noticeable properties\footnote{Recall that $g_{jk}\equiv g(e_{j})\cdot
e_{k}=G(e_{j},e_{k})\equiv G_{jk}$ and $g^{jk}\equiv g^{-1}(e^{j})\cdot
e^{k}=G^{jk}$ are the $jk$-entries of the inverse matrix for $[G_{jk}]$.}
\begin{align}
h(e_{j})\underset{\eta}{\cdot}h(e_{k}) &  =g_{jk,}\label{GB.2a}\\
h^{\ast}(e^{j})\underset{\eta}{\cdot}h^{\ast}(e^{k}) &  =g^{jk}.\label{GB.2b}%
\end{align}
The bases $\{h(e_{k})\}$ and $\{h^{\ast}(e^{k})\}$ are called the \emph{gauge
bases} associated to $\{e_{k}\}$ and $\{e^{k}\}.$

\subsection{Tetrad Bases}

Let $u_{1},\ldots,u_{n}$ be the $n$ \emph{euclidean} orthonormal eigenvectors
of $\eta$, i.e., the eigenvalues equation for $\eta$ can be written
(reordering $u_{1},\ldots,u_{n}$ if necessary) as
\[
\eta(u_{k})=\left\{
\begin{array}
[c]{cc}%
u_{k}, & k=1,\ldots,p\\
-u_{k}, & k=p+1,\ldots,n
\end{array}
\right.  ,
\]
and $u_{j}\cdot u_{k}=\delta_{jk}.$

The $h^{-1}$-deformed bases $\{h^{-1}(u_{k})\}$ and $\{h^{\dagger}(u_{k})\}$
satisfy the remarkable properties
\begin{align}
h^{-1}(u_{j})\underset{g}{\cdot}h^{-1}(u_{k}) &  =\eta_{jk},\label{TB.1a}\\
h^{\dagger}(u_{j})\underset{g^{-1}}{\cdot}h^{\dagger}(u_{k}) &  =\eta
_{jk},\label{TB.1b}%
\end{align}
where

$\eta_{jk}\equiv\eta(u_{j})\cdot u_{k}=\left\{
\begin{array}
[c]{cc}%
1, & j=k=1,\ldots,p\\
-1, & j=k=p+1,\ldots,n\\
0, & j\equiv k
\end{array}
\right.  .$ 

The bases $\{h^{-1}(u_{k})\}$ and $\{h^{\dagger}(u_{k})\}$ are called the
\emph{tetrad bases} associated to $\{u_{k}\}.$

\subsection{Some Details of the Tetrad Formalism}

Algebraically speaking the tetrad formalism of General Relativity deals at
each tangent tensor space at a given point of a manifold with the so-called
\emph{tetrad components} of vectors, tensors, etc., i.e., the contravariant
and covariant components of vectors, tensors, etc., with respect to
\emph{tetrad bases.} A tetrad basis has the remarkable property that the
tetrad components of the metric tensor all are just real constants. Let us
recall all that.

Let $\{\partial_{i}\}$ and $\{\partial^{i}\}$ be two reciprocal bases for $V,$
i.e., $\partial_{i}\cdot\partial^{j}=\delta_{i}^{j}$. The Latin indices will
be called here \textit{coordinate indices}.

Let us take another arbitrary pair of reciprocal bases for $V,$ say
$\{e_{\alpha}\}$ and $\{e^{\alpha}\},$ i.e., $e_{\alpha}\cdot e^{\beta}%
=\delta_{\alpha}^{\beta}$ (the Greek letters are used as `tetrad indices').
Associated to it we can construct another pair of reciprocal bases for $V$ by
using the gauge extensor $h,$
\begin{align}
\varepsilon_{\alpha} &  =h^{-1}(e_{\alpha})\label{TF.1a}\\
\varepsilon^{\alpha} &  =h^{\dagger}(e^{\alpha}).\label{TF.1b}%
\end{align}
They are indeed a pair of reciprocal bases, since $h^{-1}(e_{\alpha})\cdot
h^{\dagger}(e^{\beta})=\delta_{\alpha}^{\beta},$ and are called \emph{tetrad
bases}.

As we can easily see, when $\{e_{\alpha}\}$ and $\{e^{\alpha}\}$ are thought
from the viewpoint of $\{\varepsilon_{\alpha}\}$ and $\{\varepsilon^{\alpha
}\}$ as being $e_{\alpha}=h(\varepsilon_{\alpha})$ and $e^{\alpha}=h^{\ast
}(e^{\alpha}),$ the first ones are just the gauge bases associated to the
second ones.

Consider now the contravariant and covariant components of the basis vectors
$\varepsilon_{\alpha}$ with respect to the vector bases $\{\partial^{i}\}$ and
$\{g(\partial_{i})\},$ i.e.,
\begin{align}
\left.  \varepsilon_{\alpha}\right.  ^{i}  &  =\varepsilon_{\alpha}%
\cdot\partial^{i}\label{TF.2a}\\
\varepsilon_{\alpha i}  &  =\varepsilon_{\alpha}\cdot g(\partial_{i}),
\label{TF.2b}%
\end{align}
and consider also the contravariant and covariant components of the basis
vectors $\varepsilon^{\alpha}$ with respect to the vector bases $\{g^{-1}%
(\partial^{i})\}$ and $\{\partial_{i}\},$ i.e.,
\begin{align}
\varepsilon^{\alpha i}  &  \equiv\varepsilon^{\alpha}\cdot g^{-1}(\partial
^{i})\label{TF.3a}\\
\left.  \varepsilon^{\alpha}\right.  _{i}  &  \equiv\varepsilon^{\alpha}%
\cdot\partial_{i}. \label{TF.3b}%
\end{align}

These various kinds of contravariant and covariant components for the basis
vectors $\varepsilon_{\alpha}$ and $\varepsilon^{\alpha}$ will satisfy some
well-known properties which appear in books on general relativity. Indeed, if
$g_{ij}\equiv g(\partial_{i})\cdot\partial_{j}$ and $g^{ij}\equiv
g^{-1}(\partial^{i})\cdot\partial^{j},$ we have
\begin{align}
\varepsilon_{\alpha i} &  =g_{ij}\left.  \varepsilon_{\alpha}\right.
^{j}\label{TF.4a}\\
\left.  \varepsilon_{\alpha}\right.  ^{i} &  =g^{ij}\varepsilon_{\alpha
j}.\label{TF.4b}%
\end{align}

The $n^{2}+n^{2}$ real numbers $\left.  \varepsilon_{\alpha}\right.  ^{i}$ and
$\left.  \varepsilon^{\alpha}\right.  _{i}$ are the entries of inverses
matrices to each other, i.e.,
\begin{align}
\left.  \varepsilon_{\alpha}\right.  ^{i}\left.  \varepsilon^{\alpha}\right.
_{j}  &  =\delta_{j}^{i}\label{TF.5a}\\
\left.  \varepsilon_{\alpha}\right.  ^{i}\left.  \varepsilon^{\beta}\right.
_{i}  &  =\delta_{\alpha}^{\beta}. \label{TF.5b}%
\end{align}

If $g_{\alpha\beta}\equiv g(\varepsilon_{\alpha})\cdot\varepsilon_{\beta}$ and
$g^{\alpha\beta}\equiv g^{-1}(\varepsilon^{\alpha})\cdot\varepsilon^{\beta},$
and $\eta_{\alpha\beta}\equiv\eta(e_{\alpha})\cdot e_{\beta}$ and
$\eta^{\alpha\beta}\equiv\eta^{-1}(e^{\alpha})\cdot e^{\beta},$ then
\begin{align}
g_{\alpha\beta}  &  =\eta_{\alpha\beta}=\left.  \varepsilon_{\alpha}\right.
^{i}\varepsilon_{\beta i}\label{TF.6a}\\
g^{\alpha\beta}  &  =\eta^{\alpha\beta}=\varepsilon^{\alpha i}\left.
\varepsilon^{\beta}\right.  _{i},\label{TF.6b}\\
\varepsilon_{\alpha i}  &  =\eta_{\alpha\beta}\left.  \varepsilon^{\beta
}\right.  _{i}\label{TF.6c}\\
\left.  \varepsilon^{\alpha}\right.  _{i}  &  =\eta^{\alpha\beta}%
\varepsilon_{\beta i}. \label{TF.6d}%
\end{align}

Also, the $\varepsilon^{\alpha}$-contravariant and $g(\varepsilon_{\alpha}%
)$-covariant components of a vector $v,$ i.e., $v^{\alpha}\equiv
\varepsilon^{\alpha}\cdot v$ and $v_{\alpha}\equiv g(\varepsilon_{\alpha
})\cdot v,$ can be written in terms of the $\partial^{i}$-contravariant and
$g(\partial_{i})$-covariant components of $v,$ i.e., $v^{i}\equiv
v\cdot\partial^{i}$ and $v_{i}\equiv v\cdot g(\partial_{i}),$ by the following
formulas
\begin{align}
v^{\alpha} &  =\left.  \varepsilon^{\alpha}\right.  _{i}v^{i}=\varepsilon
^{\alpha i}v_{i}\label{TF.7a}\\
v_{\alpha} &  =\left.  \varepsilon_{\alpha}\right.  ^{i}v_{i}=\varepsilon
_{\alpha i}v^{i}.\label{TF.7b}%
\end{align}

As can be easily checked, e.g.,  for a covariant $2$-tensor $T,$ the
$\varepsilon_{\alpha}$-covariant components $T_{\alpha\beta}\equiv
T(\varepsilon_{\alpha},\varepsilon_{\beta})$ and the $\partial_{i}$-covariant
components $T_{ij}\equiv T(\partial_{i},\partial_{j})$ are related by
\begin{align}
T_{\alpha\beta} &  =T_{ij}\left.  \varepsilon_{\alpha}\right.  ^{i}\left.
\varepsilon_{\beta}\right.  ^{j}\label{TF.8a}\\
T_{ij} &  =T_{\alpha\beta}\left.  \varepsilon^{\alpha}\right.  _{i}\left.
\varepsilon^{\beta}\right.  _{j}.\label{TF.8b}%
\end{align}
And the $\varepsilon^{\alpha}$-contravariant components $T^{\alpha\beta}\equiv
T(\varepsilon^{\alpha},\varepsilon^{\beta})$ and the $\partial^{i}%
$-contravariant components $T^{ij}\equiv T(\partial^{i},\partial^{j})$ are
related by
\begin{align}
T^{\alpha\beta} &  =T^{ij}\left.  \varepsilon^{\alpha}\right.  _{i}\left.
\varepsilon^{\beta}\right.  _{j}\label{TF.9a}\\
T^{ij} &  =T^{\alpha\beta}\left.  \varepsilon_{\alpha}\right.  ^{i}\left.
\varepsilon_{\beta}\right.  ^{j}.\label{TF.9b}%
\end{align}

\section{Golden Formula}

Let $h$ be any gauge extensor for $g,$ i.e., $g=h^{\dagger}\circ\eta\circ h,$
where $\eta$ is a pseudo-orthogonal metric extensor with the same signature as
$g.$ Let $\underset{g}{\ast}$ mean either $\wedge$ (exterior product),
$\underset{g}{\cdot}$ ($g$-scalar product), $\underset{g}{\lrcorner}%
,\underset{g}{\llcorner}$ ($g$-contracted products) or $\underset{g}{}$
($g$-Clifford product).  And analogously for $\underset{\eta}{\ast}.$

The $g$-metric products $\underset{g}{\ast}$ and the $\eta$-metric products
are related by a remarkable formula, called in what follows the \emph{golden
formula}. For all $X,Y\in\bigwedge V$
\begin{equation}
\underline{h}(X\underset{g}{\ast}Y)=[\underline{h}(X)\underset{\eta}{\ast
}\underline{h}(Y)],\label{GF.1}%
\end{equation}
where $\underline{h}$ denotes the \textit{extended }\cite{2} of $h$.

\textbf{Proof}

By recalling the fundamental properties for the extended of a $(1,1)$%
-extensor: $\underline{t}(X\wedge Y)=\underline{t}(X)\wedge\underline{t}(Y)$
and $\underline{t}(\alpha)=\alpha,$ we have that Eq.(\ref{GF.1}) holds for the
exterior product, i.e.,
\begin{equation}
X\wedge Y=\underline{h}^{-1}[\underline{h}(X)\wedge\underline{h}(Y)]
\label{GF.2}%
\end{equation}
and, by recalling Eq.(\ref{GB.1a}), Eq.(\ref{GF.1}) holds also for the
$g$-scalar product and the $\eta$-scalar product, i.e.,
\begin{equation}
X\underset{g}{\cdot}Y=\underline{h}^{-1}[\underline{h}(X)\underset{\eta}%
{\cdot}\underline{h}(Y)]. \label{GF.3}%
\end{equation}

By using the multivector identities for an invertible operator: $\underline
{t}^{\dagger}(X)\lrcorner Y=\underline{t}^{-1}[X\lrcorner\underline{t}(Y)]$
and $X\llcorner\underline{t}^{\dagger}(Y)=\underline{t}^{-1}[\underline
{t}(X)\llcorner Y],$ and Eq.(\ref{GE.1}) we can easily prove that
Eq.(\ref{GF.1}) holds for the $g$-contracted product and the $\eta$-contracted
product, i.e.,
\begin{align}
X\underset{g}{\lrcorner}Y  &  =\underline{h}^{-1}[\underline{h}(X)\underset
{\eta}{\lrcorner}\underline{h}(Y)]\label{GF.4a}\\
X\underset{g}{\llcorner}Y  &  =\underline{h}^{-1}[\underline{h}(X)\underset
{\eta}{\llcorner}\underline{h}(Y)]. \label{GF.4b}%
\end{align}

In order to prove Eq.(\ref{GF.4a}) recall that we can write
\[
X\underset{g}{\lrcorner}Y=\underline{h^{\dagger}\circ\eta\circ h}(X)\lrcorner
Y=\underline{h}^{-1}[\underline{\eta\circ h}(X)\lrcorner\underline
{h}(Y)]=\underline{h}^{-1}[\underline{h}(X)\underset{\eta}{\lrcorner
}\underline{h}(Y)],
\]
where the definitions of $\underset{g}{\lrcorner}$ and $\underset{\eta
}{\lrcorner}$ have been used. The proof of Eq.(\ref{GF.4b}) is completely
analogous, the definitions of $\underset{g}{\llcorner}$ and $\underset{\eta
}{\llcorner}$ should be used.

In order to prove that Eq.(\ref{GF.1}) holds for the $g$-Clifford product and
the $\eta$-Clifford product, i.e.,
\begin{equation}
X\underset{g}{}Y=\underline{h}^{-1}[\underline{h}(X)\underset{\eta}%
{}\underline{h}(Y)], \label{GF.5}%
\end{equation}
we first must prove four particular cases of it.

Take $\alpha\in R$ and $X\in\Lambda V.$ By using the axioms of the $g$ and
$\eta$ Clifford products: $\alpha\underset{g}{}X=X\underset{g}{}\alpha=\alpha
X$ and $\alpha\underset{\eta}{}X=X\underset{\eta}{}\alpha=\alpha X,$ we can
write
\[
\alpha\underset{g}{}X=\alpha X=\underline{h}^{-1}[\alpha\underline
{h}(X)]=\underline{h}^{-1}[\alpha\underset{\eta}{}\underline{h}(X)],
\]
i.e.,
\begin{equation}
\alpha\underset{g}{}X=\underline{h}^{-1}[\underline{h}(\alpha)\underset{\eta
}{}\underline{h}(X)]. \label{GF.6}%
\end{equation}
Analogously, we have
\begin{equation}
X\underset{g}{}\alpha=\underline{h}^{-1}[\underline{h}(X)\underset{\eta}%
{}\underline{h}(\alpha)]. \label{GF.7}%
\end{equation}

Take $v\in V$ and $X\in\Lambda V.$ By using the axioms of the $g$ and $\eta$
Clifford products: $v\underset{g}{}X=v\underset{g}{\lrcorner}X+v\wedge X$ and
$v\underset{\eta}{}X=v\underset{\eta}{\lrcorner}X+v\wedge X,$ and
Eqs.(\ref{GF.4a}) and (\ref{GF.2}) we can write
\[
v\underset{g}{}X=v\underset{g}{\lrcorner}X+v\wedge X=\underline{h}%
^{-1}[h(v)\underset{\eta}{\lrcorner}\underline{h}(X)]+\underline{h}%
^{-1}[h(v)\wedge\underline{h}(X)],
\]
i.e.,
\begin{equation}
v\underset{g}{}X=\underline{h}^{-1}[h(v)\underset{\eta}{}\underline{h}(X)].
\label{GF.8}%
\end{equation}
From the axioms of the $g$ and $\eta$ Clifford products: $X\underset{g}%
{}v=X\underset{g}{\llcorner}v+X\wedge v$ and $X\underset{\eta}{}%
v=X\underset{\eta}{\llcorner}v+X\wedge v$, and Eqs.(\ref{GF.4b}) and
(\ref{GF.2}) we get
\begin{equation}
X\underset{g}{}v=\underline{h}^{-1}[\underline{h}(X)\underset{\eta}{}h(v)].
\label{GF.9}%
\end{equation}

Take $v_{1},v_{2},\ldots,v_{k}\in V.$ By using $k-1$ times Eq.(\ref{GF.8}) we
have indeed that
\begin{align}
v_{1}\underset{g}{}v_{2}\underset{g}{\cdots}v_{k}  &  =\underline{h}%
^{-1}[h(v_{1})\underset{\eta}{}\underline{h}(v_{2}\underset{g}{\cdots}%
v_{k})]\nonumber\\
&  =\underline{h}^{-1}[h(v_{1})\underset{\eta}{}h(v_{2})\underset{\eta}%
{\cdots}h(v_{k})],\nonumber\\
v_{1}\underset{g}{}v_{2}\underset{g}{\cdots}v_{k}  &  =\underline{h}%
^{-1}[h(v_{1})\underset{\eta}{}h(v_{2})\underset{\eta}{\cdots}h(v_{k})].
\label{GF.10}%
\end{align}

Take $v_{1},v_{2},\ldots,v_{k}\in V$ and $X\in\Lambda V.$ By using $k-1$ times
Eq.(\ref{GF.8}) and Eq.(\ref{GF.10}) we have indeed that
\begin{align}
(v_{1}\underset{g}{}v_{2}\underset{g}{\cdots}v_{k})\underset{g}{}X  &
=\underline{h}^{-1}[h(v_{1})\underset{\eta}{}\underline{h}(\underset{g}{}%
v_{2}\underset{g}{\cdots}v_{k}\underset{g}{}X)]\nonumber\\
&  =\underline{h}^{-1}[h(v_{1})\underset{\eta}{}h(v_{2})\underset{\eta}%
{\cdots}h(v_{k})\underset{\eta}{}\underline{h}(X)],\nonumber\\
(v_{1}\underset{g}{}v_{2}\underset{g}{\cdots}v_{k})\underset{g}{}X  &
=\underline{h}^{-1}[\underline{h}(v_{1}\underset{g}{}v_{2}\underset{g}{\cdots
}v_{k})\underset{\eta}{}\underline{h}(X)]. \label{GF.11}%
\end{align}

We now can prove the general case of Eq.(\ref{GF.5}). We shall use an
expansion formula for multivectors: $X=X^{0}+\underset{k=1}{\overset{n}{\sum}%
}\dfrac{1}{k!}X^{j_{1}\ldots j_{k}}e_{j_{1}}\underset{g}{\cdots}e_{j_{k}},$
where $\{e_{j}\}$ is any basis for $V,$ Eq.(\ref{GF.6}) and Eq.(\ref{GF.11}).
We can write
\begin{align*}
X\underset{g}{}Y  &  =X^{0}\underset{g}{}Y+\underset{k=1}{\overset{n}{\sum}%
}\frac{1}{k!}X^{j_{1}\ldots j_{k}}(e_{j_{1}}\underset{g}{\cdots}e_{j_{k}%
})\underset{g}{}Y\\
&  =\underline{h}^{-1}[\underline{h}(X^{0})\underset{\eta}{}\underline
{h}(Y)]+\underline{h}^{-1}[\underset{k=1}{\overset{n}{\sum}}\frac{1}%
{k!}X^{j_{1}\ldots j_{k}}\underline{h}(e_{j_{1}}\underset{g}{\cdots}e_{j_{k}%
})\underset{\eta}{}\underline{h}(Y)]\\
&  =\underline{h}^{-1}[\underline{h}(X^{0}+\underset{k=1}{\overset{n}{\sum}%
}\frac{1}{k!}X^{j_{1}\ldots j_{k}}e_{j_{1}}\underset{g}{\cdots}e_{j_{k}%
})\underset{\eta}{}\underline{h}(Y)],\\
X\underset{g}{}Y  &  =\underline{h}^{-1}[\underline{h}(X)\underset{\eta}%
{}\underline{h}(Y)].
\end{align*}

Hence, Eq.(\ref{GF.2}), Eq.(\ref{GF.3}), Eqs.(\ref{GF.4a}) and (\ref{GF.4b}),
and Eq.(\ref{GF.5}) have set the \emph{golden }formula.$\blacksquare$

It should be noticed that the relationship between the $g^{-1}$-metric
products $\underset{g^{-1}}{*}$ and the $\eta$-metric products $\underset
{\eta}{*}$ is given by
\begin{equation}
\underline{h}^{*}(X\underset{g^{-1}}{*}Y)=\underline{h}^{*}(X)\underset{\eta
}{*}\underline{h}^{*}(Y). \label{GF.12}%
\end{equation}

\section{Conclusions}

In this paper we just continued the program started at \cite{1} towards the
construction of a theory of geometrical algebra of multivectors and the theory
of extensors. In Section 2 we introduced the concepts of metric extensor and
of metric adjoint operators. In Section 3 we introduced pseudo-orthogonal
metric extensors and in particular the important case of Lorentz extensors.
Gauge extensors, gauge bases, tetrad bases and some algebraic aspects of the
tetrad formalism (fundamental tools in the formulation of geometric theories
of the gravitational field\cite{rodoliv2006}) are studied in Sections 4 and 5.
In section 6 we prove the remarkable \emph{golden formula}, which permit us to
do calculations in an arbitrary metric Clifford algebra $\mathcal{C}\ell(V,G)$
in terms of the euclidean algebra $\mathcal{C}\ell(V,G_{E})$ since the former
algebra is interpreted as a precise \textit{gauge} deformation of the later.
This idea is at the basis of our formulation of the theory of deformed
geometries, to be introduced in following papers of the series.\medskip

\textbf{Acknowledgments: }V. V. Fern\'{a}ndez and A. M. Moya are very grateful
to Mrs. Rosa I. Fern\'{a}ndez who gave to them material and spiritual support
at the starting time of their research work. This paper could not have been
written without her inestimable help.

\end{document}